\documentclass[11pt,a4paper]{article}
\usepackage[utf8]{inputenc}
\usepackage[english]{babel}
\usepackage[T1]{fontenc} 
\usepackage{lmodern} 
\usepackage{xcolor} 
\usepackage{color}
\usepackage{amssymb} 
\usepackage{mathrsfs} 
\usepackage{amsmath}
\usepackage{amsthm}
\usepackage{amsopn}
\usepackage{enumitem}
\usepackage{graphicx} 
\usepackage{url}
\usepackage{hyperref}
\usepackage{doi}
\usepackage[margin=3cm]{geometry}
\usepackage[titletoc]{appendix}

\theoremstyle{plain}

\newtheorem{theorem}{Theorem}[section]

\theoremstyle{definition}

\newtheorem{example}[theorem]{Example}

\theoremstyle{remark}

\newtheorem{remark}[theorem]{Remark}

\newcommand{\N}{{\mathbb N}}

\newcommand{\R}{{\mathbb R}}

\newcommand{\eps}{\varepsilon}
\newcommand{\e}{\mathrm{e}}
\def\1{{\mathchoice {1\mskip-4mu\mathrm l}      
{1\mskip-4mu\mathrm l}
{1\mskip-4.5mu\mathrm l} {1\mskip-5mu\mathrm l}}}
\newcommand{\measurerestr}{%
  \,\raisebox{-.127ex}{\reflectbox{\rotatebox[origin=br]{-90}{$\lnot$}}}\,%
}
\DeclareMathOperator{\restrict}{\hspace{-0.2ex}\measurerestr\hspace{-0.2ex}}

\newcommand{\weak}{\stackrel{\scriptstyle w}\longrightarrow}
\newcommand{\weakstar}{\stackrel{\ast}\rightharpoonup}

\begin{document}
\begin{center}
\begin{Large}
Most likely balls in Banach spaces: existence and non-existence
\end{Large}
\\[0.5cm]
\begin{large}
Bernd Schmidt\footnote{Universität Augsburg, Germany, {\tt bernd.schmidt@math.uni-augsburg.de}} 
\end{large}
\\[0.5cm]
\today
\\[1cm]
\end{center}

\begin{abstract}
We establish a general criterion for the existence of convex sets of fixed shape as, e.g., balls of a given radius, of maximal probability on Banach spaces. We also provide counterexamples showing that their existence my fail even in some common situations.  
\end{abstract}

\medskip 


\section{Introduction}

In this note we address the natural question if for a (Borel) probability measure $\mu$ on a separable Banach space $X$ and for a given radius there is always a ``ball of maximum likelihood'', i.e., if the maximum of $\mu(B)$ is attained among all balls $B \subset X$ of radius $r > 0$. 
For small radii, such a maximizer---if existent---can be viewed as an approximation to and a regularization of a ``mode'' referring to a ``point of maximum likelihood'' for the given measure $\mu$ as $r$ becomes small. 
More generally, instead of balls of a fixed radius our main existence theorem will also apply to the system of all translates of a fixed convex set $C$. 

This issue has received considerable attention recently, notably in the area of Bayesian inverse problems (cf.\ \cite{Stuart:10}) where the problem arises to devise maximum a posteriori (MAP) estimators, in particular for Gaussian priors on infinite dimensional (separable) Banach spaces, \cite{DashtiLawStuartVoss:13,HelinBurger:15,LieSullivan:18,ClasonHelinKretschmannPiiroinen:19,KlebanovWacker:23,Lambley:23}. 
Indeed, seminal results in that area made implicit use of the existence of balls of maximum likelihood (cp.\ the discussion in \cite{KlebanovWacker:23}), while very recently it has been noted that the question of their existence can be circumvented by considering the asymptotics of almost maximizers in order to obtain MAP estimators \cite{KlebanovWacker:23,Lambley:23}. 
Nevertheless, besides being a question of intrinsic interest, the problem remains relevant as the quest for the position of a (small) ball of positive radius with maximal probability amounts to solving a regularized point optimization problem. 
As such it typically enjoys improved stability properties and might even be favorable from a modeling perspective. 
In particular this will be the case in situations when, possibly due to data uncertainties in the presence of noise, it is preferable to estimate parameter regions rather than single points. 

The problem is addressed for general metric spaces in some detail in \cite[Sect.~4.2]{LambleySullivan:22}. 
There the authors provide a collection of technical sufficient conditions for the existence of maximum likelihood balls (``radius-$r$ modes'' in their terminology). 
They also give counterexamples for some particular measures on specific spaces and certain ranges of $r$. 
The farthest reaching existence result for measures on Banach spaces known to date seems to be \cite[Thm.~4.9]{LambleySullivan:22}, which proves that on the sequence spaces $\ell^p$, $1 < p < \infty$, maximum likelihood balls do exist for measures that do not charge any sphere in $\ell^p$. 
In particular this applies to ``radius-$r$ maximum likelihood a posteriori estimators for Gaussian priors'' on $\ell^p$, $1 < p < \infty$, which are absolutely continuous with respect to a non-degenerate Gaussian. 
It appears that no counterexample on a Banach space is known to date. 

The purpose of this note is twofold. 
First, we establish a general existence theorem for maximum likelihood convex shapes (and, in particular, balls of any radius) on Banach spaces. 
In particular this will apply to every separable and reflexive space (and $\ell^1$), thus closing a gap in the seminal contribution \cite{DashtiLawStuartVoss:13}.  
Second, by way of various examples we also show that existence my fail in some natural situations. 
In fact, we will provide a couple of counterexamples on $c_0$ and the Wiener space, which might even be absolutely continuous with respect to a non-degenerate Gaussian measure, and which do not allow for balls of maximal probability for any value of radius $r$.

\section{A general existence result}

Throughout we assume that $X$ is a separable (real) Banach space. 
By $B_r$ and $B_r^{\circ}$ we denote the closed and, respectively, open ball of radius $r$ in $X$ centered at $0$.  
For $x \in X$, $C \subset X$ we write $x + C =: C(x)$. 
Suppose $\mu$ is a Borel probability measure on $X$. 
   
\begin{theorem}
Suppose $X$ is the separable dual of a Banach space. Let $C \subset X$ be a bounded weak*-closed convex set $C \subset X$ with non-empty interior. Then there exists $x_0 \in X$ such that 
$$ \mu(C(x_0)) \ge \mu(C(x)) $$ 
for all $x \in X$. 
\end{theorem}

\begin{remark}
\begin{enumerate}
 \item In particular, this applies to every separable reflexive Banach space $X$. 
 \item An admissible choice for $C$ is $C = B_r$ for any $r > 0$. 
\end{enumerate}
\end{remark}

\begin{proof}
Without loss of generality we assume $0$ is an interior point of $C$. 
We consider a maximizing sequence of translates $C(x_n)$, i.e., 
\begin{align*} 
 \mu(C(x_n)) \to m_0 := \sup \big\{\mu(C(x)) : x \in X \big\}. 
\end{align*}
Since $X$ is separable and $C$ has non-empty interior, it is easy to see that $m_0 > 0$. 
Clearly, the sequence $(x_n)$ is bounded. (If $R > 0$ is chosen such that $\mu(X \setminus B_R) \le \frac{m_0}{2}$, then $x_n \in B_{R + \operatorname{diam}(C)}$ eventually.) As $\mu$ is tight (on the polish space $X$), also the family $(\mu_n)$ of restrictions $\mu_n = \mu \restrict C(x_n)$ is tight. Thus Prohorov's theorem implies $\mu_n \weak \mu_0$ weakly (in duality with $C_b(X)$) for a (not relabeled) subsequence and a finite measure $\mu_0$. It follows that 
$$ \mu_0(X) 
   = \lim_{n \to \infty} \mu_n(X) = \lim_{n \to \infty} \mu(C(x_n)) = m_0. $$ 

By assumption, $X$ has a separable predual $X_*$, so by Alaoglu's theorem we may pass to a further subsequence (not relabeled) such that $x_n \weakstar x_0$ weakly* in $X$ for some $x_0 \in X$. 
We will now prove that $\mu_0$ is supported on $C(x_0)$. 

To this end, we fix any $z \notin C(x_0)$. Since $C$ is weak*-closed, with the help of the Hahn-Banach theorem for the dual pairing $(X, X_*)$ we can choose an element $x_* \in X_*$ and then an $\eps > 0$ such that
\begin{align}\label{HB}
  \sup \big\{ \langle x_*, y \rangle : y \in B_\eps(z) \big\} 
  = \langle x_*, z \rangle + \eps \| x_* \| 
  < \inf \big\{ \langle x_*, y \rangle : y \in C(x_0) \big\}. 
\end{align}
On the other hand, the Portmanteau theorem implies 
$$ \mu_0(B^{\circ}_\eps(z))
   \le \liminf_{n \to \infty} \mu_n(B^{\circ}_\eps(z))
   = \liminf_{n \to \infty} \mu(C(x_n) \cap B^{\circ}_\eps(z)). $$ 
As a consequence we conclude that, in case $\mu_0(B^{\circ}_\eps(z)) > 0$ we have $C(x_n) \cap B^{\circ}_\eps(z)) \ne \emptyset$ for sufficiently large $n$, say $x_n + y_n \in C(x_n) \cap B^{\circ}_\eps(z))$ (and so $y_n\in C$). Passing to yet another subsequence (not relabeled) we get $y_n \weakstar y$ for some $y \in C$. It follows that $x_n + y_n \weakstar x_0 + y \in C(x_0) \cap B_\eps(z)$, which contradicts \eqref{HB}. So we must have $\mu_0(B^{\circ}_\eps(z)) = 0$. This proves that $\operatorname{supp} \mu_0 \subset C(x_0)$. 

It remains to observe that $\mu_0 \le \mu$, which follows from the outer regularity of the Borel measures $\mu_0$ and $\mu$ and from the fact that for any open subset $U \subset X$ the Portmanteau theorem gives  
$$ \mu_0(U) 
   \le \liminf_{n \to \infty} \mu_n(U) 
   = \liminf_{n \to \infty} \mu(C(x_n) \cap U) 
   \le \mu(U). $$ 
Summarizing we find that 
$$ \mu(C(x_0)) 
   \ge \mu_0(C(x_0)) 
   = \mu_0(X) 
   \ge m_0, $$
which, by definition of $m_0$, proves $\mu(C(x_0)) = m_0$. 
\end{proof}

\section{Examples of non-existence}

We discuss a number of concrete cases, where balls of maximum likelihood do not exist. 
There is a common underlying idea in all of them which would easily allow to generate further examples along these lines. 

Our first two examples are on the Banach space $X = c_0$ of (real) null-sequences equipped with the $\sup$-norm, which is separable and even has a separable dual (namely, $\ell^1$), but is not a dual space itself. 

\begin{example}\label{ex:Exp}
We let $\mu = \bigotimes_{k \in \N} \mathrm{Exp}(k)$, where $\mathrm{Exp}(\lambda)$ denotes the exponential distribution on $\R$ with rate parameter $\lambda$ (and cumulative distribution function $x \mapsto 1 - \e^{-\lambda x^+}$). An easy application of the Borel-Cantelli lemma shows $\mu(X) = 1$. Let $r > 0$ arbitrary. 

For any ball $B_r(x)$, $x = (x_1, x_2, \ldots)$ one has 
$$ \mu(B_r(x)) 
   = \prod_{k\in\N} \big( \e^{-k(x_k-r)^+} - \e^{-k(x_k+r)^+} \big). $$ 
Choosing $k_0$ such that $x_{k_0} < r$ and setting $x'_k = x_k$ for $k \ne k_0$, $x'_{k_0} = r$ we get $\mu(B_r(x')) > \mu(B_r(x))$ if $\mu(B_r(x)) > 0$. This shows that $x \mapsto \mu(B_r(x))$ does not have a maximizer. (Its supremum is $m_0 = \prod_{k\in\N} \big( 1 - \e^{-2kr} \big)$, which can be seen by maximizing each factor separately and considering the maximizing sequence $(x_{\cdot,n}) \subset c_0$, $x_{k,n} = r$ for $k \le n$, $x_{k,n} = 0$ for $k > n$.) 
\end{example}

\begin{example}
In order to give an example where $\mu$ is absolutely continuous with respect to a Gaussian measure on $X = c_0$, we first let $\mu_0 = \bigotimes_{k \in \N} \mathcal{N}(0, k^{-2})$, where $\mathcal{N}(0, \sigma^2)$ denotes the Gaussian on $\R$ with variance mean $0$ and $\sigma^2$ and cumulative distribution function $X \mapsto \Phi(x/\sigma)$. An easy application of the Borel-Cantelli lemma shows $\mu_0(X) = 1$. We consider the (closed) set $A \subset X$ given by 
$$ A = \big\{x \in X : x_k \ge - 1/\sqrt{k} \text{ for all } k \in \N\big\} $$ 
and note that $\mu_0(A)$ (the probability that the coordinate process does not pass the moving boundary $k \mapsto -1/\sqrt{k}$) is positive since $ \mu_0(X \setminus A) \le \sum_{k \in \N} \Phi(-\sqrt{k}) < 1$.
We then define $\mu$ by conditioning on $A$, i.e., we set $\mu = \frac{1}{\mu(A)} \1_A \mu_0$. 

A similar reasoning as above shows that balls of maximal probability do not exist (and the supremum is explicitly given as $m_0 = \prod_{k\in\N} [ \Phi((kr - \sqrt{k})^++kr) - \Phi((kr - \sqrt{k})^+-kr) ] $).
\end{example}

We now give some examples on the Wiener space $X = \{\omega \in C[0,1] : \omega(0) = 0\}$ equipped, as usual, with the $\sup$-norm in order to show that non-existence of maximum likelihood balls is encountered in common situations in a continuous time setting. They follow the similar basic idea of the previous two examples. 

In what follows we let $\mu_0$ be the Wiener measure on $X$ so that the coordinate process $(\omega(t))_{t \in [0,1]}$ is Brownian motion. 

\begin{example}
Similarly as in Example~\ref{ex:Exp} we can consider typical processes which assume only non-negative values as, e.g., the running maximum of Brownian motion or the reflected Brownian motion:  
$$ \omega_{\max}(t) = \max\{\omega(s) : 0 \le s \le t \},
\quad \text{respectively,} \quad |\omega|(t) = |\omega(t)|. 
$$ 
If $\mu$ denotes the corresponding distribution on $X$, in both cases the maximum of $x \mapsto \mu(B_r(x))$ is not attained. Indeed, as $\omega(t) \to 0$ for $t \to 0$ for any $\omega \in X$ and so $\omega(t) < r$ on $[0,s]$ for some $0 < s < 1$, it suffices to choose any $\omega' \in X$ such that $\omega < \omega' \le r$ on $(0,s)$ and $\omega' = \omega$ on $[t, 1]$ to get $\mu(B_r(\omega')) > \mu(B_r(\omega))$ if $\mu(B_r(\omega)) > 0$.   
\end{example}

\begin{example}
We choose a non-positive $\rho \in X$ such that $\mu_0(\{\omega : \omega(t) \ge \rho(t) \text{ for all }t \in [0,1]\}) > 0$. Such a $\rho$ can be found with the help of Khinchin's law of the iterated logarithm: $\liminf_{t\to 0} \omega(t)/\sqrt{2t \log \log (1/t)} = -1$ for $\mu_0$-a.e.\ $\omega \in X$, which allows to choose $0 < t_0 < 1/\e$ such that 
$$ \mu_0\Big(\Big\{\omega(t) \ge -2 \sqrt{2t \log \log (1/t)}\Big\}\Big) \text{ for all } t \in [0,t_0] 
   > 0. $$
We now define $\rho \in X$ by $\rho(t) = -2\sqrt{2t \log \log (1/t)}$ for $t \le t_0$ and then $\rho(t) = \rho(t_0)$ for $t > t_0$. The (closed) set $A \subset X$ 
$$ A = \big\{\omega \in X : \omega(t) \ge \rho(t)  \text{ for all } t \in [0,1]\big\} \subset X 
$$ 
will then have positive probability $\mu_0(A) > 0$. 
Conditioning on $A$, we define $\mu = \frac{1}{\mu(A)} \1_A \mu_0$. 

For any ball $B_r(\omega)$ we have $\mu(B_r(\omega)) > 0$ if and only if $\omega + r > \rho$ on $[0,1]$. A construction as in the previous example shows that $\mu(B_r(\omega')) > \mu(B_r(\omega))$ for such $\omega$.  So, again, the supremum of these values is not attained. 
\end{example}

\subsubsection*{Acknowledgments} 
I am grateful to Philipp Wacker for drawing my attention to this problem and interesting conversations on the subject. 

\bibliographystyle{alpha} 
\bibliography{Literature}
\end{document}